\documentclass{amsart}
\usepackage{amssymb}
\usepackage{amsfonts}

\newtheorem{theorem}{Theorem}
\theoremstyle{plain}

\newtheorem{corollary}{Corollary}

\newtheorem{definition}{Definition}

\newtheorem{lemma}{Lemma}

\newtheorem{proposition}{Proposition}
\newtheorem{remark}{Remark}

\numberwithin{equation}{section}
\input{tcilatex}

\begin{document}
\title[\textbf{TAYLOR SPECTRUM}]{\textbf{TAYLOR SPECTRUM AND CHARACTERISTIC
FUNCTIONS OF COMMUTING \ 2-CONTRACTIONS}}
\author{BERRABAH BENDOUKHA}
\email{bbendoukha@gmail.com}
\address{University of mostaganem\\
Department of Mathematics\\
B.O. 227, Mostaganem (27000), Algeria}

\begin{abstract}
In this paper, we give a description of Taylor spectrum of commuting
2-contractions in terms of characteritic functions of such contractions. The
case of a single contraction obtained by B. Sz. Nagy and C. Foias is
generalied in this work.
\end{abstract}

\maketitle

\section{Introduction}

Let $\mathcal{H}$ be a Hilbert space. A $2$-tuple $A=\left(
A_{1},A_{2}\right) $ of bounded operators on $\mathcal{H}$ is called
contractive ( or $2$-contraction ) if $\left\Vert
A_{1}h_{1}+A_{2}h_{2}\right\Vert ^{2}\leq \left\Vert h_{1}\right\Vert
^{2}+\left\Vert h_{2}\right\Vert ^{2}$ for all $h_{1},h_{2}$ in $\mathcal{H}$%
. It is equivalent \cite{1} to the condition: $A_{1}A_{1}^{\star
}+A_{2}A_{2}^{\star }\leq 1_{\mathcal{H}}$. If additionaly, operators $A_{1}$
and $A_{2}$ commute, then $A$ is called a commuting $2$-contraction. To
every commuting 2-contraction $A=\left( A_{1},A_{2}\right) $ corresponds an
analytic operator-valued function $\ \theta _{A}:ID^{2}\rightarrow B\left( 
\mathcal{D}_{A}\text{ , }\mathcal{D}_{A^{\ast }}\right) $ called
characteristic function of $A$ and defined by : \ 
\begin{equation}
\theta _{A}\left( z_{1},z_{2}\right) =-A+D_{A^{\ast }}\left( 1_{\mathcal{H}%
}-z_{1}A_{1}^{\ast }-z_{2}A_{2}^{\ast }\right) ^{-1}\left( z_{1}.1_{\mathcal{%
H}},z_{2}.1_{\mathcal{H}}\right) D_{A}  \label{eq1}
\end{equation}%
\newline
where

\begin{description}
\item[(a)] $\mathbb{D}^{2}=\left\{ \left( z_{1},z_{2}\right) \in 
\mathbb{C}
^{2}:\left\vert z_{1}\right\vert ^{2}+\left\vert z_{1}\right\vert
^{2}<1\right\} ,$

\item[(b)] $B\left( \mathcal{D}_{A}\text{ , }\mathcal{D}_{A^{\ast }}\right) $
is the set of all bounded operators from $\mathcal{D}_{A}$ into $\mathcal{D}%
_{A^{\ast }},$

\item[(c)] $D_{A^{\ast }}=\left( I-A_{1}A_{1}^{\star }-A_{2}A_{2}^{\star
}\right) ^{\frac{1}{2}}:\mathcal{H}\rightarrow \mathcal{H}$ and $\mathcal{D}%
_{A^{\ast }}$ is the closure of the range of $D_{A^{\ast }}$,\newline

\item[(d)] $D_{A}=\left[ 
\begin{array}{cc}
1_{H}-A_{1}^{\ast }A_{1} & -A_{1}^{\ast }A_{2} \\ 
-A_{2}^{\ast }A_{1} & 1_{H}-A_{2}^{\ast }A_{2}%
\end{array}%
\right] ^{\frac{1}{2}}:\mathcal{H}^{2}\rightarrow \mathcal{H}^{2}$ and is
the closure of the range of $D_{A}$,\newline

\item[(e)] $\left( z_{1}1_{\mathcal{H}},z_{2}1_{\mathcal{H}}\right) :%
\mathcal{H}^{2}\rightarrow \mathcal{H};\qquad \left( z_{1}1_{\mathcal{H}%
},z_{2}1_{\mathcal{H}}\right) \left( 
\begin{array}{c}
h_{1} \\ 
h_{2}%
\end{array}%
\right) =z_{1}h_{1}+z_{2}h_{2}$.\newline
\end{description}

Characteristic function of commuting $n$-contraction has been introduced in 
\cite{3} as a generalization of characteristic function of a single
contraction \cite{9}. A lot of its remarquable properties have been
established in (\cite{2}, \cite{4},\cite{9}). In particular, it is shown
(like in the single case) that the characteristic function is a unitary
invariant. It means that characteristic functions of two pure or completely
noncoisometric $n$-contractions (\cite{2}, \cite{4}) $A=\left(
A_{1},...,A_{n}\right) $ and $A^{\prime }=\left( A_{1}^{\prime
},...,A_{n}^{\prime }\right) $ coincide if and only if there exists a
unitary operator $U:\mathcal{H}\rightarrow \mathcal{H}$ such that $%
A_{i}^{\prime }=U^{-1}A_{i}U$ \ for every $i=1,...,n$.

In the case of a completly nonunitary single contraction, the spectrum can
be described in terms of the characteristic function (see theoerm 4.1 \cite%
{9}). The aim of this paper is to give a description of Taylor spectrum in
the case of commuting pure 2-contractions by means of characteristic
function (\ref{eq1}).

In section 2 we briefly remind the definition of Taylor spectrum. Section 3
contains\ characterizations of different components of Taylor spectrum. In
section 4, we investigate the behavior of Taylor spectrum under the action
of involutive automorphims of unit ball.

\section{Taylor spectrum}

Let $A=\left( A_{1},A_{2},...,A_{n}\right) $ be a pure $n$-contraction.
According to (\cite{5}, \cite{6}, \cite{7}, \cite{10}), the Taylor spectrum
of $A$ can be defined as follows. Let $\Lambda \left( \mathcal{H}\right) $
be the exterior algebra on $n$ generators $e_{1}$,..., $e_{n}$ with identity 
$e_{0}=1$ and coefficients in $H$. In other words, 
\begin{equation*}
\Lambda \left( \mathcal{H}\right) =\left\{ x\otimes e_{i_{1}}\wedge
...\wedge e_{i_{p}}:x\in \mathcal{H};\,1\leq i_{1}....\lessdot i_{p}\leq n%
\text{ };\text{ }1\leq p\leq n\right\}
\end{equation*}%
with the collapsing property : $e_{i}\wedge e_{j}+e_{j}\wedge e_{i}=0$. One
has 
\begin{eqnarray*}
\Lambda \left( \mathcal{H}\right) &=&\oplus _{k=1}^{n}\Lambda ^{k}\left( 
\mathcal{H}\right) \text{; \qquad }\Lambda ^{k}\left( \mathcal{H}\right)
=\left\{ x\otimes e_{i_{1}}\wedge ...\wedge e_{i_{k}}:x\in H,\,1\leq
i_{1}....\lessdot i_{k}\leq n\right\} \text{; } \\
\Lambda ^{0}\left( \mathcal{H}\right) &=&\mathcal{H}.
\end{eqnarray*}%
Consider in $\Lambda \left( \mathcal{H}\right) $ operator: 
\begin{equation*}
B_{A}:\Lambda \left( \mathcal{H}\right) \rightarrow \Lambda \left( \mathcal{H%
}\right) :\text{\qquad }B_{A}\left( x\otimes e_{i_{1}}\wedge ...\wedge
e_{i_{p}}\right) =\dsum\limits_{k=1}^{n}A_{k}\left( x\right) \otimes
e_{k}\wedge e_{i_{1}}\wedge ...\wedge e_{i_{p}}.
\end{equation*}%
\newline
It is not difficult to see that $B_{A}^{2}=0$ and $Ran\,B_{A}\subseteq
Ker\,B_{A}$ . Decomposition $\Lambda \left( \mathcal{H}\right) =\oplus
_{k=1}^{n}\Lambda ^{k}\left( \mathcal{H}\right) $ gives a rise to a cochain $%
K\left( A,\mathcal{H}\right) $, the so-called Koszul complex $K\left( A,%
\mathcal{H}\right) $ associated to $A$ on $\mathcal{H}$ as follows: 
\begin{equation*}
K\left( A_{1},A_{2},\mathcal{H}\right) :\left\{ 0\right\} \rightarrow 
\mathcal{H}\,=\Lambda ^{0}\left( \mathcal{H}\right) \,^{\underrightarrow{%
B_{A}^{0}}}\,\,\,\,\,....\,\,\,^{^{\underrightarrow{B_{A}^{n-1}}%
}\,\,}\Lambda ^{n}\left( \mathcal{H}\right) \,\rightarrow \left\{ 0\right\}
\end{equation*}%
where $B_{A}^{k}$ \ is the restriction of $B_{A}$ to the subspace $\Lambda
^{k}\left( \mathcal{H}\right) $. Complex $K\left( A,\mathcal{H}\right) $ is
said to be exact (or regular) if: 
\begin{equation*}
\left\{ 0\right\} =\ker \,B_{A}^{0}\text{, \ }Ran\,B_{A}^{0}=\ker \,B_{A}^{1}%
\text{,......, }Ran\,B_{A}^{n-2}=\ker \,B_{A}^{n-1}\text{, \ }%
Ran\,B_{A}^{n-1}=\Lambda ^{n}\left( H\right) .
\end{equation*}

\begin{definition}
\label{def1}The Taylor spectrum of$\ n$-contraction is the set: 
\begin{equation*}
\sigma _{T}\left( A\right) =\left\{ z=\left(
z_{1},\,z_{2,....,\,}z_{n}\right) \in 
\mathbb{C}
^{n}\text{ : }K\left( A_{1}-z_{1},....,A_{n}-z_{n};H\right) \text{ is not
exact}\right\} .
\end{equation*}%
Let us now suppose that $n=2$. Then, 
\begin{eqnarray*}
\Lambda \left( \mathcal{H}\right) &=&\Lambda ^{0}\left( \mathcal{H}\right)
\oplus \Lambda ^{1}\left( \mathcal{H}\right) \oplus \Lambda ^{2}\left( 
\mathcal{H}\right) \\
&=&\left( \mathcal{H}\otimes e_{0}\right) \oplus \left( \left( \mathcal{H}%
\otimes e_{1}\right) \oplus \left( \mathcal{H}\otimes e_{2}\right) \right)
\oplus \left( \mathcal{H}\otimes e_{1}\wedge e_{2}\right) .
\end{eqnarray*}
\end{definition}

According to this direct sum, operator $B_{A}$ admits the matrix
representation 
\begin{equation*}
B_{A}=\left[ 
\begin{array}{cccc}
0 & 0 & 0 & 0 \\ 
A_{1} & 0 & 0 & 0 \\ 
A_{2} & 0 & 0 & 0 \\ 
0 & -A_{2} & A_{1} & 0%
\end{array}%
\right]
\end{equation*}%
and then operators $B_{\left( A_{1},A_{2}\right) }^{0}$ and $B_{\left(
A_{1},A_{2}\right) }^{1}$ have the forms:%
\begin{equation}
\left\{ 
\begin{array}{c}
B_{\left( A_{1},A_{2}\right) }^{0}\left( x\right) =A_{1}\left( x\right)
\oplus A_{2}\left( x\right) ,\text{ \ \ \ }\left( x,y\in \mathcal{H}\right) ,
\\ 
B_{\left( A_{1},A_{2}\right) }^{1}\left( x\oplus y\right) =-A_{2}\left(
x\right) +A_{1}\left( y\right) ,\text{ \ \ }\left( x,y\in \mathcal{H}\right)%
\end{array}%
\right. .  \label{eq2}
\end{equation}

According to definition\ \ref{def1} and formula (\ref{eq2}), one has 
\begin{equation*}
\sigma _{T}\left( A_{1},A_{2}\right) =\sigma _{T}^{\left( 1\right) }\left(
A_{1},A_{2}\right) \cup \sigma _{T}^{\left( 2\right) }\left(
A_{1},A_{2}\right) \cup \sigma _{T}^{\left( 3\right) }\left(
A_{1},A_{2}\right)
\end{equation*}%
where%
\begin{equation}
\left( z_{1},z_{2}\right) \in \sigma _{T}^{\left( 1\right) }\left(
A_{1},A_{2}\right) \Leftrightarrow \exists x\in \mathcal{H}:x\neq 0,\ \
\left( A_{1}-z_{1}\right) x=\left( A_{2}-z_{2}\right) x=0  \label{eq3}
\end{equation}%
\begin{equation}
\left( z_{1},z_{2}\right) \in \sigma _{T}^{\left( 2\right) }\left(
A_{1},A_{2}\right) \Leftrightarrow \exists \left( x_{1},x_{2}\right) \in 
\mathcal{H}^{2}:\left\{ 
\begin{array}{c}
\left( A_{1}-z_{1}\right) x_{1}-\left( A_{2}-z_{2}\right) x_{2}=0 \\ 
\left( x_{1},x_{2}\right) \neq \left( \left( A_{1}-z_{1}\right) h\text{ },%
\text{ }\left( A_{2}-z_{2}\right) h\right) \text{, }\forall h\in \mathcal{H}%
\end{array}%
\right.  \label{eq4}
\end{equation}%
\begin{equation}
\left( z_{1},z_{2}\right) \in \sigma _{T}^{\left( 3\right) }\left(
A_{1},A_{2}\right) \Leftrightarrow \exists y\in \mathcal{H}:y\neq \left(
A_{1}-z_{1}\right) x_{1}-\left( A_{2}-z_{2}\right) x_{2},\ \ \forall \left(
x_{1},x_{2}\right) \in \mathcal{H}^{2}.  \label{eq5}
\end{equation}

\begin{remark}
$\sigma _{T}^{\left( 1\right) }\left( A_{1},A_{2}\right) $ is called the
ponctual joint Taylor spectrum. \newline
Taylor joint spectrum generalizes the one variable notion of spectrum. It is
a nonempty compact subset of $%
\mathbb{C}
^{n}$. The reader can find an excellent account of the Taylor spectrum and
its relations with other multiparameter spectral theories in \cite{6}. Note
also that in \cite{2} a description of Harte spectrum by means of
characteristic function is given.
\end{remark}

Throughout this paper, we will suppose that if $A=\left( A_{1},A_{2}\right) $
is a commuting 2-contraction, then operator $D_{A^{\ast }}=\left(
I-A_{1}A_{1}^{\star }-A_{2}A_{2}^{\star }\right) ^{\frac{1}{2}}$ is one to
one. Note that pure ( and more generally completely non coisometric )
commutig 2-contractions (\cite{2}, \cite{3}, \cite{4}) satisfy this
condition. Indeed, if $A=\left( A_{1},A_{2}\right) $ is a pure 2-contraction
then, the decreasing sequence of positive bounded operators $\left( \left(
A_{1}A_{1}^{\star }+A_{2}A_{2}^{\star }\right) ^{n}\right) _{n\in IN}$ \
admits a strong limit $A_{\infty }=0$. Because of that,

\begin{eqnarray*}
D_{A^{\ast }}\left( x\right) &=&0\Rightarrow D_{A^{\ast }}^{2}\left(
x\right) =\left( I-A_{1}A_{1}^{\star }-A_{2}A_{2}^{\star }\right) x=0 \\
&\Rightarrow &x=\left( A_{1}A_{1}^{\star }+A_{2}A_{2}^{\star }\right)
x\Rightarrow x=\left( A_{1}A_{1}^{\star }+A_{2}A_{2}^{\star }\right) ^{n}x%
\text{ , }\forall n=0,1,2,... \\
&\Rightarrow &x=A_{\infty }\left( x\right) =0\text{.}
\end{eqnarray*}

Using relations $AD_{A}=D_{A^{\ast }}A$ \ and \ $A^{\ast }D_{A^{\ast
}}=D_{A}A^{\ast }$ (\cite{3}), it can be proven that :

$D_{A^{\ast }}=\left( I-A_{1}A_{1}^{\star }-A_{2}A_{2}^{\star }\right) ^{%
\frac{1}{2}}$ is one to one $\Leftrightarrow $ $D_{A}=\left[ 
\begin{array}{cc}
1_{H}-A_{1}^{\ast }A_{1} & -A_{1}^{\ast }A_{2} \\ 
-A_{2}^{\ast }A_{1} & 1_{H}-A_{2}^{\ast }A_{2}%
\end{array}%
\right] ^{\frac{1}{2}}$ is one to one.

\section{\textbf{Characterization of Taylor spectrum}}

\begin{lemma}
\label{lem1}Let $A=\left( A_{1},A_{2}\right) $ be a commuting $2$%
-contraction such that $D_{A^{\ast }}$ is one to one. Then, 
\begin{equation*}
A\left( x,y\right) =z_{1}x+z_{2}y,\left( \left( z_{1},z_{2}\right) \in 
\mathbb{C}
^{2}\text{, }\left( x,y\right) \in \mathcal{H}^{2}\right) \Leftrightarrow
\theta _{A}\left( z_{1},z_{2}\right) \left( D_{A}\left( 
\begin{array}{c}
x \\ 
y%
\end{array}%
\right) \right) =0.
\end{equation*}
\end{lemma}

\begin{proof}
It follows directly from relation 
\begin{equation}
\theta _{A}\left( z_{1},z_{2}\right) D_{A}\left( 
\begin{array}{c}
x \\ 
y%
\end{array}%
\right) =D_{A^{\ast }}\left( 1_{\mathcal{H}}-z_{1}A_{1}^{\ast
}-z_{2}A_{2}^{\ast }\right) ^{-1}\left[ \left( z_{1}x+z_{2}y\right) -\left(
A_{1}x+A_{2}y\right) \right]  \label{eq6}
\end{equation}
\end{proof}

\begin{lemma}
\label{lem2}Let $A=\left( A_{1},A_{2}\right) $ be a commuting pure $2$%
-contraction such that $D_{A^{\ast }}$ is one to one, $\left(
z_{1},z_{2}\right) \in 
\mathbb{C}
^{2}$ and $x\in \mathcal{H}$. Then, 
\begin{equation*}
A^{\ast }\left( x\right) =\left( 
\begin{array}{c}
\overline{z_{1}}.x \\ 
\overline{z_{2}}.x%
\end{array}%
\right) \Leftrightarrow \left( \theta _{A}\left( z_{1},z_{2}\right) \right)
^{\ast }D_{A^{\ast }}\left( x\right) =0.
\end{equation*}
\end{lemma}

\begin{proof}
Since in this case $D_{A^{\ast }}^{2}\left( x\right) =\left( 1_{\mathcal{H}}-%
\overline{z_{1}}A_{1}-\overline{z_{2}}A_{2}\right) x$, then the necessary
condition is a direct conequence of relation, \ $\qquad $%
\begin{equation}
\left( \theta _{A}\left( z_{1},z_{2}\right) \right) ^{\ast }D_{A^{\ast
}}\left( x\right) =D_{A}\left( -\left( 
\begin{array}{c}
A_{1}^{\ast }x \\ 
A_{2}^{\ast }x%
\end{array}%
\right) +\left( 
\begin{array}{c}
\overline{z_{1}} \\ 
\overline{z_{2}}%
\end{array}%
\right) \left( \left( 1_{H}-\overline{z_{1}}A_{1}-\overline{z_{2}}%
A_{2}\right) ^{-1}\right) D_{A^{\ast }}^{2}\left( x\right) \right)
\label{eq7}
\end{equation}%
\newline
\underline{\textbf{Proof of the suffisant condition}} \ Hence $\left( \theta
_{A}\left( 0,0\right) \right) ^{\ast }D_{A^{\ast }}\left( x\right) =-\left( 
\begin{array}{c}
A_{1}^{\ast }x \\ 
A_{2}^{\ast }x%
\end{array}%
\right) $, one can whithout loosing the generality suppose that $\left(
z_{1},z_{2}\right) \neq \left( 0,0\right) $.%
\begin{eqnarray*}
\left( \theta _{A}\left( z_{1},z_{2}\right) \right) ^{\ast }D_{A^{\ast
}}\left( x\right) &=&0 \\
&\Rightarrow &D_{A}\left( -\left( 
\begin{array}{c}
A_{1}^{\ast }x \\ 
A_{2}^{\ast }x%
\end{array}%
\right) +\left( 
\begin{array}{c}
\overline{z_{1}} \\ 
\overline{z_{2}}%
\end{array}%
\right) \left( \left( 1_{H}-\overline{z_{1}}A_{1}-\overline{z_{2}}%
A_{2}\right) ^{-1}\right) D_{A^{\ast }}^{2}\left( x\right) \right) =0 \\
&\Rightarrow &\left\{ 
\begin{array}{c}
-A_{1}^{\ast }x+\overline{z_{1}}\left( \left( 1_{H}-\overline{z_{1}}A_{1}-%
\overline{z_{2}}A_{2}\right) ^{-1}\right) \left( I-A_{1}A_{1}^{\ast
}-A_{2}A_{2}^{\ast }\right) x=0 \\ 
-A_{2}^{\ast }x+\overline{z_{2}}\left( \left( 1_{H}-\overline{z_{1}}A_{1}-%
\overline{z_{2}}A_{2}\right) ^{-1}\right) \left( I-A_{1}A_{1}^{\ast
}-A_{2}A_{2}^{\ast }\right) x=0%
\end{array}%
\right. \\
&\Rightarrow &\left\{ 
\begin{array}{c}
-A_{1}^{\ast }x+\overline{z_{1}}\left( \left( 1_{H}-\overline{z_{1}}A_{1}-%
\overline{z_{2}}A_{2}\right) ^{-1}\right) \left( I-A_{1}A_{1}^{\ast
}-A_{2}A_{2}^{\ast }\right) x=0 \\ 
-A_{2}^{\ast }x+\overline{z_{2}}\left( \left( 1_{H}-\overline{z_{1}}A_{1}-%
\overline{z_{2}}A_{2}\right) ^{-1}\right) \left( I-A_{1}A_{1}^{\ast
}-A_{2}A_{2}^{\ast }\right) x=0%
\end{array}%
\right. \\
&\Rightarrow &\left\{ 
\begin{array}{c}
-\overline{z_{2}}A_{1}^{\ast }x+\overline{z_{2}}\overline{z_{1}}\left(
\left( 1_{H}-\overline{z_{1}}A_{1}-\overline{z_{2}}A_{2}\right) ^{-1}\right)
\left( I-A_{1}A_{1}^{\ast }-A_{2}A_{2}^{\ast }\right) x=0 \\ 
-\overline{z_{1}}A_{2}^{\ast }x+\overline{z_{1}}\overline{z_{2}}\left(
\left( 1_{H}-\overline{z_{1}}A_{1}-\overline{z_{2}}A_{2}\right) ^{-1}\right)
\left( I-A_{1}A_{1}^{\ast }-A_{2}A_{2}^{\ast }\right) x=0%
\end{array}%
\right. \\
&\Rightarrow &\left\{ 
\begin{array}{c}
\overline{z_{1}}A_{2}^{\ast }x=\overline{z_{2}}A_{1}^{\ast }x \\ 
-\overline{z_{2}}A_{1}^{\ast }x+\overline{z_{2}}\overline{z_{1}}\left(
\left( 1_{H}-\overline{z_{1}}A_{1}-\overline{z_{2}}A_{2}\right) ^{-1}\right)
\left( I-A_{1}A_{1}^{\ast }-A_{2}A_{2}^{\ast }\right) x=0%
\end{array}%
\right. \\
&\Rightarrow &\left\{ 
\begin{array}{c}
\overline{z_{1}}A_{2}^{\ast }x=\overline{z_{2}}A_{1}^{\ast }x \\ 
-\overline{z_{2}}A_{1}^{\ast }x+\overline{z_{2}}\left( \left( 1_{H}-%
\overline{z_{1}}A_{1}-\overline{z_{2}}A_{2}\right) ^{-1}\right) \left( 
\overline{z_{1}}I-\overline{z_{1}}A_{1}A_{1}^{\ast }-\overline{z_{2}}%
A_{2}A_{1}^{\ast }\right) x=0%
\end{array}%
\right. \\
&\Rightarrow &\left\{ 
\begin{array}{c}
\overline{z_{1}}A_{2}^{\ast }x=\overline{z_{2}}A_{1}^{\ast }x \\ 
\overline{z_{2}}\left( \left( 1_{H}-\overline{z_{1}}A_{1}-\overline{z_{2}}%
A_{2}\right) ^{-1}\right) \left( \left( \overline{z_{1}}I-A_{1}^{\ast
}\right) \right) x=0%
\end{array}%
\right. \\
&\Rightarrow &\left( \left( \overline{z_{1}}I-A_{1}^{\ast }\right) \right)
x=0\Rightarrow A_{1}^{\ast }x=\overline{z_{1}}x.
\end{eqnarray*}

On the other hand,

\begin{enumerate}
\item $\ \overline{z_{1}}A_{2}^{\ast }x=\overline{z_{2}}A_{1}^{\ast }x$, $\
A_{1}^{\ast }x=\overline{z_{1}}x$ \ and \ $z_{1}\neq 0\Rightarrow
A_{2}^{\ast }x=\overline{z_{2}}x.$

\item Putting $A_{1}^{\ast }x=\overline{z_{1}}x$ \ and $z_{1}=0$ in the
relation 
\begin{equation*}
-A_{2}^{\ast }x+\overline{z_{2}}\left( \left( 1_{H}-\overline{z_{1}}A_{1}-%
\overline{z_{2}}A_{2}\right) ^{-1}\right) \left( I-A_{1}A_{1}^{\ast
}-A_{2}A_{2}^{\ast }\right) x=0,
\end{equation*}%
and multiplying by $\left( 1_{H}-\overline{z_{2}}A_{2}\right) ^{-1}$, one
obtains 
\begin{equation*}
\left( 1_{H}-\overline{z_{2}}A_{2}\right) A_{2}^{\ast }x=\overline{z_{2}}%
\left( I-A_{2}A_{2}^{\ast }\right) x.
\end{equation*}%
This last relation is equivalent to $A_{2}^{\ast }x=\overline{z_{2}}x$ .
\end{enumerate}
\end{proof}

\begin{proposition}
Let $A=\left( A_{1},A_{2}\right) $ be a commuting pure $2$-contraction such
that $D_{A^{\ast }}$ is one to one. Then, $\left( z_{1},z_{2}\right) \in
\sigma _{T}^{\left( 1\right) }\left( A_{1},A_{2}\right) $ if and only if
equation $\theta _{A}\left( z_{1},z_{2}\right) X=0$ admits at least two
nontrivial solutions $D_{A}\left( X_{1}\right) $ and $D_{A}\left(
X_{2}\right) $ such that, $X_{1}=\left( 
\begin{array}{c}
x \\ 
0%
\end{array}%
\right) $ , $\ X_{2}=\left( 
\begin{array}{c}
0 \\ 
x%
\end{array}%
\right) $, $x\in \mathcal{H}$.
\end{proposition}

\begin{proof}
One has 
\begin{eqnarray*}
\left( z_{1},z_{2}\right) &\in &\sigma _{T}^{\left( 1\right) }\left(
A_{1},A_{2}\right) \\
&\Leftrightarrow &\exists x\in \mathcal{H}:x\neq 0,\left( A_{1}-z_{1}\right)
x=\left( A_{2}-z_{2}\right) x=0 \\
&\Leftrightarrow &\exists x\in \mathcal{H}:x\neq 0,\text{ }A_{1}\left(
x\right) =z_{1}.x\text{ \ and }A_{2}\left( x\right) =z_{2}.x \\
&\Leftrightarrow &\exists x\in \mathcal{H}:x\neq 0,\text{ }A\left( 
\begin{array}{c}
x \\ 
0%
\end{array}%
\right) =z_{1}.x+z_{2}.0\text{ \ and \ }A\left( 
\begin{array}{c}
0 \\ 
x%
\end{array}%
\right) =z_{1}.0+z_{2}.x \\
&\Leftrightarrow &\theta _{A}\left( z_{1},z_{2}\right) \left( D_{A}\left( 
\begin{array}{c}
x \\ 
0%
\end{array}%
\right) \right) =0\text{ \ \ and \ }\theta _{A}\left( z_{1},z_{2}\right)
\left( D_{A}\left( 
\begin{array}{c}
0 \\ 
x%
\end{array}%
\right) \right) =0
\end{eqnarray*}%
\newline
To end the proof, it is sufficient to remark that 
\begin{equation*}
x\neq 0\Leftrightarrow D_{A}\left( 
\begin{array}{c}
x \\ 
0%
\end{array}%
\right) \neq 0\Leftrightarrow D_{A}\left( 
\begin{array}{c}
0 \\ 
x%
\end{array}%
\right) \neq 0.
\end{equation*}
\end{proof}

\begin{proposition}
\label{prop3}Let $A=\left( A_{1},A_{2}\right) $ be a commuting pure $2$%
-contraction such that $D_{A^{\ast }}$ is one to one. Then, $\left(
z_{1},z_{2}\right) \in \sigma _{T}^{\left( 2\right) }\left(
A_{1},A_{2}\right) $ if and only if the equation $\theta _{A}\left(
z_{1},z_{2}\right) Y=0$ admits at least one non trivial solution $Y=D_{A}$ $%
\left( X\right) $ such that $X\neq \left( 
\begin{array}{c}
\left( A_{1}-z_{1}\right) h \\ 
\left( z_{2}-A_{2}\right) h%
\end{array}%
\right) ,$ $\forall h\in \mathcal{H}$.
\end{proposition}

\begin{proof}
One has,%
\begin{equation*}
\left( z_{1},z_{2}\right) \in \sigma _{T}^{\left( 2\right) }\left(
A_{1},A_{2}\right) \Leftrightarrow \left\{ 
\begin{array}{c}
\exists \left( x_{1},x_{2}\right) \in \mathcal{H}^{2}:\text{ \ }\left(
A_{1}-z_{1}\right) x_{1}-\left( A_{2}-z_{2}\right) x_{2}=0, \\ 
\left( x_{1},x_{2}\right) \neq \left( \left( A_{1}-z_{1}\right) h\text{ },%
\text{ }\left( A_{2}-z_{2}\right) h\right) \text{ }\forall h\in \mathcal{H},
\\ 
\left( x_{1},x_{2}\right) \neq \left( 0,0\right) .%
\end{array}%
\right.
\end{equation*}%
\newline
It means that\newline
\begin{equation*}
\left( z_{1},z_{2}\right) \in \sigma _{T}^{\left( 2\right) }\left(
A_{1},A_{2}\right) \Leftrightarrow \left\{ 
\begin{array}{c}
\exists \left( x_{1},x_{2}\right) \in \mathcal{H}^{2}:\text{ \ }A_{1}\left(
x_{1}\right) +A_{2}\left( -x_{2}\right) =z_{1}.x_{1}+z_{2}.\left(
-x_{2}\right) , \\ 
\left( x_{1},x_{2}\right) \neq \left( \left( A_{1}-z_{1}\right) h\text{ },%
\text{ }\left( A_{2}-z_{2}\right) h\right) ;\text{ }\forall h\in \mathcal{H},
\\ 
\left( x_{1},x_{2}\right) \neq \left( 0,0\right) .%
\end{array}%
\right.
\end{equation*}%
According to lemma \ref{lem1}, one has finally%
\begin{equation*}
\left( z_{1},z_{2}\right) \in \sigma _{T}^{\left( 2\right) }\left(
A_{1},A_{2}\right) \Leftrightarrow \left\{ 
\begin{array}{c}
\exists \left( x_{1},x_{2}\right) \in \mathcal{H}^{2}:\theta _{A}\left(
z_{1},z_{2}\right) \left( D_{A}\left( 
\begin{array}{c}
x_{1} \\ 
-x_{2}%
\end{array}%
\right) \right) =0, \\ 
\left( x_{1},x_{2}\right) \neq \left( \left( A_{1}-z_{1}\right) h\text{ },%
\text{ }\left( A_{2}-z_{2}\right) h\right) ;\text{ }\forall h\in \mathcal{H},
\\ 
\left( x_{1},x_{2}\right) \neq \left( 0,0\right) .%
\end{array}%
\right. .
\end{equation*}
\end{proof}

\begin{proposition}
\label{prop2}Let $A=\left( A_{1},A_{2}\right) $ be a commuting $2$%
-contraction such that $D_{A^{\ast }}$ is one to one and \textbf{\ }$\left(
z_{1},z_{2}\right) \in 
\mathbb{C}
^{2}$. Assume that equation $\left( \theta _{A}\left( z_{1},z_{2}\right)
\right) ^{\ast }D_{A^{\ast }}\left( y\right) =0$ admits at least one non
trivial solution. Then, 
\begin{equation*}
\left( z_{1},z_{2}\right) \in \sigma _{T}^{\left( 3\right) }\left(
A_{1},A_{2}\right) .
\end{equation*}
\end{proposition}

\begin{proof}
Suppose that $\left( \theta _{A}\left( z_{1},z_{2}\right) \right) ^{\ast
}D_{A^{\ast }}\left( y\right) =0$ admits at least one non trivial solution $%
y $. According to lemma \ref{lem2}, it means that 
\begin{equation*}
A^{\ast }\left( y\right) =\left( 
\begin{array}{c}
\overline{z_{1}}.y \\ 
\overline{z_{2}}.y%
\end{array}%
\right)
\end{equation*}%
\newline
Thus, for every $\left( x_{1},x_{2}\right) \in \mathcal{H}^{2}$, one has 
\begin{eqnarray*}
0 &=&\left\langle A^{\ast }\left( y\right) -\left( 
\begin{array}{c}
\overline{z_{1}}.y \\ 
\overline{z_{2}}.y%
\end{array}%
\right) ,\left( 
\begin{array}{c}
x_{2} \\ 
-x_{1}%
\end{array}%
\right) \right\rangle _{\mathcal{H}^{2}} \\
&=&\left\langle \left( 
\begin{array}{c}
\left( A_{1}^{\ast }-\overline{z_{1}}\right) .y \\ 
\left( A_{2}^{\ast }-\overline{z_{2}}\right) .y%
\end{array}%
\right) ,\left( 
\begin{array}{c}
x_{2} \\ 
-x_{1}%
\end{array}%
\right) \right\rangle _{\mathcal{H}^{2}} \\
&=&\left\langle \left( A_{1}^{\ast }-\overline{z_{1}}\right)
.y,x_{2}\right\rangle _{\mathcal{H}}+\left\langle \left( A_{2}^{\ast }-%
\overline{z_{2}}\right) .y,-x_{1}\right\rangle _{\mathcal{H}} \\
&=&\left\langle y,\left( A_{1}-z_{1}\right) x_{2}\right\rangle _{\mathcal{H}%
}-\left\langle .y,\left( A_{2}-z_{2}\right) x_{1}\right\rangle _{\mathcal{H}}
\\
&=&\left\langle y,\left( A_{1}-z_{1}\right) x_{2}-\left( A_{2}-z_{2}\right)
x_{1}\right\rangle _{\mathcal{H}}.
\end{eqnarray*}%
\newline
Thus,\ 
\begin{equation*}
y\notin RanD_{\left( A_{1}-z_{1},A_{2}-z_{2}\right) }^{1}
\end{equation*}%
and finally 
\begin{equation*}
\left( z_{1},z_{2}\right) \in \sigma _{T}^{\left( 3\right) }\left(
A_{1},A_{2}\right) .\newline
\end{equation*}
\end{proof}

\section{Taylor spectrum and involutive automorphisms of unit ball}

In this section we invetigate the Taylor spectrum under action of involutive
automorphisms of unit ball $\mathbb{D}^{2}$. Such automorphisms are defined
in \cite{8} by 
\begin{equation*}
\Phi _{\lambda }\left( z\right) =\lambda -\frac{\sqrt{1-\left\Vert \lambda
\right\Vert ^{2}}}{1-\langle z,\lambda \rangle }\left( z-\left( 1-\sqrt{%
1-\left\Vert \lambda \right\Vert ^{2}}\right) \frac{\langle z,\lambda
\rangle }{\left\Vert \lambda \right\Vert ^{2}}.\lambda \right) \text{, }%
\lambda \in \mathbb{D}^{2}\text{, }\lambda =\left( \lambda _{1},\lambda
_{2}\right) \neq 0.
\end{equation*}

Connection between automorphisms of unit ball and multicontractions has been
made in \cite{2} and \cite{4} where some very interesting properties have
been established. In particular, if $\ \Phi _{\lambda }$ is an involutive
automorphism of unit ball and $A=\left( A_{1},A_{2}\right) $ is a
commutative $2$-contraction then, (see \cite{2}, sections 4 and 5) one can
define operator 
\begin{equation*}
\Phi _{\lambda }\left( A\right) =\Lambda -D_{\Lambda ^{\ast }}\left(
1_{H}-A\Lambda ^{\ast }\right) ^{-1}AD_{\Lambda }
\end{equation*}%
where the operator $\Lambda =\left( \lambda _{1}.1_{H},\lambda
_{2}.1_{H}\right) $ is defined from $H^{2}$ into $\mathcal{H}$ by : 
\begin{equation*}
\Lambda \left( x_{1},x_{2}\right) =\lambda _{1}.x_{1}+\lambda _{2}.x_{2}.%
\newline
\end{equation*}%
Propsition \ref{prop4.1} and theorem 4.2 below summarise important for us
results obtained in \cite{2} and \cite{4}.\newline

\begin{proposition}
\label{prop4.1}Let $\Phi _{\lambda }$ an involutive automorphism of unit
ball and $A=\left( A_{1},A_{2}\right) $ a commutative 2-contraction. Then, $%
\Phi _{\lambda }\left( A\right) $ is a commutative 2-contraction such that,%
\begin{eqnarray}
\text{ \ \ \ }I-\Phi _{\lambda }\left( A\right) ^{\ast }\Phi _{\lambda
}\left( A\right) &=&D_{\Lambda }\left( 1_{\mathcal{H}}-A^{\ast }\Lambda
\right) ^{-1}\left( I-A^{\ast }A\right) \left( 1_{\mathcal{H}}-\Lambda
^{\ast }A\right) ^{-1}D_{\Lambda },  \label{eq8} \\
I-\Phi _{\lambda }\left( A\right) \Phi _{\lambda }\left( A\right) ^{\ast }
&=&D_{\Lambda ^{\ast }}\left( 1_{\mathcal{H}}-A\Lambda ^{\ast }\right)
^{-1}\left( I-AA^{\ast }\right) \left( 1_{\mathcal{H}}-\Lambda A^{\ast
}\right) ^{-1}D_{\Lambda ^{\ast }}\newline
.  \label{eq9}
\end{eqnarray}
\end{proposition}

\begin{theorem}
\label{th1}Let $\Phi _{\lambda }$ be an involutive automorphism of unit ball
and $A=\left( A_{1},A_{2}\right) $ a commutative 2-contraction. Then,

\begin{enumerate}
\item Operators \ $\Omega :$\DH $_{\Phi _{\lambda }\left( A\right)
}\rightarrow $\DH $_{A}$ and $\Omega _{\ast }:$\DH $_{\Phi _{\lambda }\left(
A\right) ^{\ast }}\rightarrow $\DH $_{A^{\ast }}$ defined by 
\begin{eqnarray*}
\Omega \left( D_{\Phi _{\lambda }\left( A\right) }\left( X\right) \right)
&=&D_{A}\left( 1_{\mathcal{H}}-\Lambda ^{\ast }A\right) ^{-1}D_{\Lambda
}\left( X\right) ,\text{ } \\
\Omega _{\ast }\left( D_{\Phi _{\lambda }\left( A\right) ^{\ast }}\left(
X\right) \right) &=&D_{A^{\ast }}\left( 1_{\mathcal{H}}-\Lambda A^{\ast
}\right) ^{-1}D_{\Lambda ^{\ast }}\left( X\right)
\end{eqnarray*}%
\newline
are unitaries.

\item $\theta _{\Phi _{\lambda }\left( A\right) }$ and $\theta _{A}$ are
connected by the relation 
\begin{equation*}
\Omega _{\ast }\theta _{\Phi _{\lambda }\left( A\right) }\left(
z_{1},z_{2}\right) =\theta _{A}\left( \Phi _{\lambda }\left(
z_{1},z_{2}\right) \right) \Omega .\newline
\end{equation*}
\end{enumerate}
\end{theorem}

It can be shown that:\newline
\begin{equation}
\left( 1_{\mathcal{H}}-A\Lambda ^{\ast }\right) ^{-1}=\left( 1_{\mathcal{H}}-%
\overline{z_{1}}A_{1}-\overline{z_{2}}A_{2}\right) ^{-1},  \label{eq10}
\end{equation}%
\begin{equation}
\left( 1_{\mathcal{H}^{2}}-\Lambda ^{\ast }A\right) ^{-1}=\left[ 
\begin{array}{cc}
\left( 1_{\mathcal{H}}-\overline{\lambda _{1}}A_{1}-\overline{\lambda _{2}}%
A_{2}\right) ^{-1}\left( 1_{\mathcal{H}}-\overline{\lambda _{2}}A_{2}\right)
& \overline{\lambda _{2}}A_{1}.\left( 1_{\mathcal{H}}-\overline{\lambda _{1}}%
A_{1}-\overline{\lambda _{2}}A_{2}\right) ^{-1} \\ 
\overline{\lambda _{1}}A_{2}.\left( 1_{\mathcal{H}}-\overline{\lambda _{1}}%
A_{1}-\overline{\lambda _{2}}A_{2}\right) ^{-1} & \left( 1_{\mathcal{H}}-%
\overline{\lambda _{1}}A_{1}-\overline{\lambda _{2}}A_{2}\right) ^{-1}\left(
1_{\mathcal{H}}-\overline{\lambda _{1}}A_{1}\right)%
\end{array}%
\right]  \label{eq11}
\end{equation}%
\begin{equation}
D_{\Lambda }=\frac{1}{\left\vert \lambda _{1}\right\vert ^{2}+\left\vert
\lambda _{2}\right\vert ^{2}}\left[ 
\begin{array}{cc}
\left\vert \lambda _{2}\right\vert ^{2}+\left\vert \lambda _{1}\right\vert
^{2}\sqrt{1-\left\vert \lambda _{2}\right\vert ^{2}-\left\vert \lambda
_{1}\right\vert ^{2}} & \overline{\lambda _{1}}\lambda _{2}\left( \sqrt{%
1-\left\vert \lambda _{2}\right\vert ^{2}-\left\vert \lambda _{1}\right\vert
^{2}}-1\right) \\ 
\lambda _{1}\overline{\lambda _{2}}\left( \sqrt{1-\left\vert \lambda
_{2}\right\vert ^{2}-\left\vert \lambda _{1}\right\vert ^{2}}-1\right) & 
\left\vert \lambda _{1}\right\vert ^{2}+\left\vert \lambda _{2}\right\vert
^{2}\sqrt{1-\left\vert \lambda _{2}\right\vert ^{2}-\left\vert \lambda
_{1}\right\vert ^{2}}%
\end{array}%
\right] .  \label{eq12}
\end{equation}

Using (\ref{eq10}), one can show that 
\begin{equation}
\Phi _{\lambda }\left( A\right) =\left( B_{1}\left( \lambda \right)
,B_{2}\left( \lambda \right) \right)  \label{eq13}
\end{equation}%
where%
\begin{eqnarray*}
B_{1}\left( \lambda \right) &=&\lambda _{1}.1_{\mathcal{H}}-\sqrt{%
1-\left\vert \lambda _{2}\right\vert ^{2}-\left\vert \lambda _{1}\right\vert
^{2}}\left( 1_{H}-\overline{\lambda _{1}}A_{1}-\overline{\lambda _{2}}%
A_{2}\right) ^{-1} \\
&&\times \left\{ A_{1}\left( \left\vert \lambda _{2}\right\vert
^{2}+\left\vert \lambda _{1}\right\vert ^{2}\sqrt{1-\left\vert \lambda
_{2}\right\vert ^{2}-\left\vert \lambda _{1}\right\vert ^{2}}\right)
+\lambda _{1}\overline{\lambda _{2}}A_{2}\left( \sqrt{1-\left\vert \lambda
_{2}\right\vert ^{2}-\left\vert \lambda _{1}\right\vert ^{2}}-1\right)
\right\}
\end{eqnarray*}%
\begin{eqnarray*}
B_{2}\left( \lambda \right) &=&\lambda _{2}.1_{\mathcal{H}}-\sqrt{%
1-\left\vert \lambda _{2}\right\vert ^{2}-\left\vert \lambda _{1}\right\vert
^{2}}\left( 1_{H}-\overline{\lambda _{1}}A_{1}-\overline{\lambda _{2}}%
A_{2}\right) ^{-1} \\
&&\times \left\{ \overline{\lambda _{1}}\lambda _{2}A_{1}\left( \sqrt{%
1-\left\vert \lambda _{2}\right\vert ^{2}-\left\vert \lambda _{1}\right\vert
^{2}}-1\right) +A_{2}\left( \left\vert \lambda _{1}\right\vert
^{2}+\left\vert \lambda _{2}\right\vert ^{2}\sqrt{1-\left\vert \lambda
_{2}\right\vert ^{2}-\left\vert \lambda _{1}\right\vert ^{2}}\right)
\right\} .
\end{eqnarray*}

Formulas (\ref{eq10}), (\ref{eq11}), and (\ref{eq12}) allow us to find the
explicit forms of operators $\Phi _{\lambda }\left( A\right) $, $\Omega $
and $\Omega _{\ast }$. On the other hand, from (\ref{eq9})\ follows that if $%
D_{A^{\ast }}$ is one to one, then\ $D_{\Phi _{\lambda }\left( A\right)
^{\ast }}$ is also one to one. Using theorem\ \ref{th1}, one can obtain the
following caracterization for Taylor spectrum of $\Phi _{\lambda }\left(
A\right) $ in terms of solutions of equations 
\begin{equation*}
\theta _{A}\left( z_{1},z_{2}\right) D_{A}\left( X\right) =0\text{ \ \ \ and
\ \ \ }\left( \theta _{A}\left( z_{1},z_{2}\right) \right) ^{\ast
}D_{A^{\ast }}\left( y\right) =0.\newline
\end{equation*}

\begin{proposition}
$\Phi _{\lambda }\left( z_{1},z_{2}\right) \in \sigma _{T}^{\left( 1\right)
}\left( \Phi _{\lambda }\left( A\right) \right) $ if and only if equation $%
\theta _{A}\left( z_{1},z_{2}\right) D_{A}\left( X\right) =0$ admits at
least two nontrivial solutions $X_{1}$ and $X_{2}$ such that, 
\begin{equation*}
X_{1}=\left( 1_{\mathcal{H}}-\Lambda ^{\ast }A\right) ^{-1}D_{\Lambda
}\left( 
\begin{array}{c}
y \\ 
0%
\end{array}%
\right) ,\text{\ \ }X_{2}=\left( 1_{\mathcal{H}}-\Lambda ^{\ast }A\right)
^{-1}D_{\Lambda }\left( 
\begin{array}{c}
0 \\ 
y%
\end{array}%
\right) ,\text{ \ }y\in \mathcal{H}.
\end{equation*}
\end{proposition}

\begin{proof}
According to Proposition \ref{prop2}, $\Phi _{\lambda }\left(
z_{1},z_{2}\right) \in \sigma _{T}^{\left( 1\right) }\left( \Phi _{\lambda
}\left( A\right) \right) $ if and only if there exists a nonnul vector $y\in 
\mathcal{H}$ such that, 
\begin{equation*}
Y_{1}=\left( 
\begin{array}{c}
y \\ 
0%
\end{array}%
\right) ,\text{ \ }Y_{2}=\left( 
\begin{array}{c}
0 \\ 
y%
\end{array}%
\right)
\end{equation*}%
are solutions of equation 
\begin{equation*}
\theta _{\Phi _{\lambda }\left( A\right) }\left( \Phi _{\lambda }\left(
z_{1},z_{2}\right) \right) D_{\Phi _{\lambda }\left( A\right) }\left(
Y\right) =0.\newline
\end{equation*}%
Using theorema \ref{th1} and the fact that $\Phi _{\lambda }$ is involutive,
it is equivalent to the existence of a nonnul vector $y\in \mathcal{H}$ such
that, 
\begin{equation*}
Y_{1}=\left( 
\begin{array}{c}
y \\ 
0%
\end{array}%
\right) ,\text{ \ }Y_{2}=\left( 
\begin{array}{c}
0 \\ 
y%
\end{array}%
\right)
\end{equation*}%
are solutions of \ equation 
\begin{equation*}
\Omega _{\ast }^{-1}\theta _{A}\left( z_{1},z_{2}\right) \Omega D_{\Phi
_{\lambda }\left( A\right) }\left( Y\right) =\Omega _{\ast }^{-1}\theta
_{A}\left( z_{1},z_{2}\right) D_{A}\left( \left( 1_{\mathcal{H}}-\Lambda
^{\ast }A\right) ^{-1}D_{\Lambda }\left( Y\right) \right) =0
\end{equation*}%
which is equivalent to the equation \newline
\begin{equation*}
\theta _{A}\left( z_{1},z_{2}\right) D_{A}\left( \left( 1_{\mathcal{H}%
}-\Lambda ^{\ast }A\right) ^{-1}D_{\Lambda }\left( Y\right) \right) =0.
\end{equation*}
\end{proof}

\begin{corollary}
Let $\left( z_{1},z_{2}\right) \in \mathbb{C}^{2}$. Then, $\left(
z_{1},z_{2}\right) \in $ $\sigma _{T}^{\left( 1\right) }\left( A\right) $
and $\Phi _{\lambda }\left( z_{1},z_{2}\right) \in $ $\sigma _{T}^{\left(
1\right) }\left( \Phi _{\lambda }\left( A\right) \right) $ if and only if
there exists two nonnul vectors $x$ and $y$ in $\mathcal{H}$ such that
vectors: 
\begin{equation*}
X_{1}=\left( 
\begin{array}{c}
x \\ 
0%
\end{array}%
\right) ,\ \text{\ \ }X_{2}=\left( 
\begin{array}{c}
0 \\ 
x%
\end{array}%
\right) ,
\end{equation*}%
\begin{equation*}
Y_{1}=\left( 1_{\mathcal{H}}-\Lambda ^{\ast }A\right) ^{-1}D_{\Lambda
}\left( 
\begin{array}{c}
y \\ 
0%
\end{array}%
\right) \text{ \ \ and \ \ \ }Y_{2}=\left( 1_{\mathcal{H}}-\Lambda ^{\ast
}A\right) ^{-1}D_{\Lambda }\left( 
\begin{array}{c}
0 \\ 
y%
\end{array}%
\right)
\end{equation*}%
are both solutions of equation 
\begin{equation*}
\theta _{A}\left( z_{1},z_{2}\right) D_{A}\left( X\right) =0.\newline
\end{equation*}
\end{corollary}

\begin{proposition}
Let\textbf{\ }$\left( z_{1},z_{2}\right) \in \mathbb{C}^{2}$. Assume that
equation $\left( \theta _{A}\left( z_{1},z_{2}\right) \right) ^{\ast
}D_{A^{\ast }}\left( y\right) =0$ admits at least one non trivial solution.
Then, 
\begin{equation*}
\Phi _{\lambda }\left( z_{1},z_{2}\right) \in \sigma _{T}^{\left( 3\right)
}\left( \Phi _{\lambda }\left( A\right) \right) .\newline
\end{equation*}
\end{proposition}

\begin{proof}
Note at first that operator $D_{\Lambda ^{\ast }}$ is invertible. Since $%
\Phi _{\lambda }$ is involutive then, according theorem \ref{th1}, 
\begin{eqnarray*}
\left( \theta _{A}\left( z_{1},z_{2}\right) \right) ^{\ast }D_{A^{\ast
}}\left( y\right) &=&0 \\
&\Leftrightarrow &\left( \theta _{A}\left( \Phi _{\lambda }\left( \Phi
_{\lambda }\left( z_{1},z_{2}\right) \right) \right) \right) ^{\ast
}D_{A^{\ast }}\left( y\right) =0 \\
&\Leftrightarrow &\left( \Omega _{\ast }\theta _{\Phi _{\lambda }\left(
A\right) }\left( \Phi _{\lambda }\left( z_{1},z_{2}\right) \right) \Omega
^{-1}\right) ^{\ast }D_{A^{\ast }}\left( y\right) =0 \\
&\Leftrightarrow &\Omega \left( \theta _{\Phi _{\lambda }\left( A\right)
}\Phi _{\lambda }\left( z_{1},z_{2}\right) \right) ^{\ast }\Omega _{\ast
}^{-1}D_{A^{\ast }}\left( y\right) =0 \\
&\Leftrightarrow &\left( \theta _{\Phi _{\lambda }\left( A\right) }\Phi
_{\lambda }\left( z_{1},z_{2}\right) \right) ^{\ast }\Omega _{\ast
}^{-1}D_{A^{\ast }}\left( y\right) =0 \\
&\Leftrightarrow &\left( \theta _{\Phi _{\lambda }\left( A\right) }\Phi
_{\lambda }\left( z_{1},z_{2}\right) \right) ^{\ast }\Omega _{\ast
}^{-1}D_{A^{\ast }}\left( \left( 1_{\mathcal{H}}-\Lambda A^{\ast }\right)
^{-1}D_{\Lambda ^{\ast }}D_{\Lambda ^{\ast }}^{-1}\left( 1_{\mathcal{H}%
}-\Lambda A^{\ast }\right) y\right) =0 \\
&\Leftrightarrow &\left( \theta _{\Phi _{\lambda }\left( A\right) }\Phi
_{\lambda }\left( z_{1},z_{2}\right) \right) ^{\ast }D_{\Phi _{\lambda
}\left( A\right) ^{\ast }}\left( D_{\Lambda ^{\ast }}^{-1}\left( 1_{\mathcal{%
H}}-\Lambda A^{\ast }\right) y\right) =0 \\
&\Leftrightarrow &\left( \theta _{\Phi _{\lambda }\left( A\right) }\Phi
_{\lambda }\left( z_{1},z_{2}\right) \right) ^{\ast }D_{\Phi _{\lambda
}\left( A\right) ^{\ast }}\left( X\right) =0
\end{eqnarray*}%
where%
\begin{equation*}
\text{ \ }X=D_{\Lambda ^{\ast }}^{-1}\left( 1_{\mathcal{H}}-\Lambda A^{\ast
}\right) y.
\end{equation*}%
Since $y$ is nonnul then, $X=D_{\Lambda ^{\ast }}^{-1}\left( 1_{H}-\Lambda
A^{\ast }\right) y$ is also nonnul and according proposition \ref{prop2}, it
follows that 
\begin{equation*}
\Phi _{\lambda }\left( z_{1},z_{2}\right) \in \sigma _{T}^{\left( 3\right)
}\left( \Phi _{\lambda }\left( A_{1},A_{2}\right) \right) .
\end{equation*}
\end{proof}

\begin{proposition}
$\Phi _{\lambda }\left( z_{1},z_{2}\right) \in \sigma _{T}^{\left( 2\right)
}\left( A_{1},A_{2}\right) $ if and only if the equation $\theta _{A}\left(
z_{1},z_{2}\right) D_{A}\left( X\right) =0$ admits at least one solution $X$
\ such that 
\begin{equation*}
X\neq \left( 
\begin{array}{c}
B_{1}\left( \lambda \right) h-w_{1}.h \\ 
B_{2}\left( \lambda \right) h-w_{2}.h%
\end{array}%
\right) ,\text{ \ \ }\forall h\in \mathcal{H}
\end{equation*}%
where $\left( w_{1},w_{2}\right) =\Phi _{\lambda }\left( z_{1},z_{2}\right) $%
.\newline
\end{proposition}

\begin{proof}
It follows immeditely from proposition \ref{prop3}.\newline
\end{proof}

\end{document}